\documentclass[reqno,onecolumn,12pt]{amsart}
\usepackage{amsfonts}
\usepackage{amssymb}
\usepackage{amsmath}
\usepackage{graphicx}
\usepackage{amscd}

\newtheorem{theorem}{Theorem}
\theoremstyle{plain}

\newtheorem{proposition}{Proposition}

\numberwithin{equation}{section}

\input{tcilatex}

\begin{document}
\author{}
\title{}
\maketitle

\begin{center}
\thispagestyle{empty} \textbf{A NEW APPROACH TO MULTIVARIATE }$q$-\textbf{%
EULER POLYNOMIALS BY USING UMBRAL CALCULUS}

\bigskip \hspace{-0.5cm}

$^{\dag }$\textbf{Serkan ARACI, }$^{\ddag }$\textbf{Xiangxing KONG, }$^{\dag
}$\textbf{Mehmet ACIKGOZ, and }$^{\sharp }$\textbf{Erdo\u{g}an \c{S}EN}

\hspace{-0.5cm}

$^{\dag }$University of Gaziantep, Faculty of Arts and Science, Department
of Mathematics, 27310 Gaziantep, Turkey

$^{\ddag }$Department of Mathematics and Statistics, Central South
University, Changsha 410075, China

$^{\sharp }$Department of Mathematics, Faculty of Science and Letters, Namik
Kemal University, 59030 Tekirda\u{g}, TURKEY

\hspace{-0.5cm}

saraci88@yahoo.com.tr; xiangxingkong@gmail.com; acikgoz@gantep.edu.tr;
erdogan.math@gmail.com

\hspace{-0.5cm}

\textbf{{\Large {Abstract}}}

\hspace{-0.5cm}
\end{center}

In the present paper, we derive numerous identities for multivariate $q$%
-Euler polynomials by using umbral calculus.

\bigskip


\hspace{-0.5cm}


\section{\textbf{Preliminaries}}

Throughout this paper, we use the following notation, where $%
\mathbb{C}
$ denotes the set of complex numbers, $\mathcal{F}$ denotes the set of all
formal power series in the variable $t$ over $%
\mathbb{C}
$ with $\mathcal{F}=\left\{ f\left( t\right) =\sum_{k=0}^{\infty }a_{k}\frac{%
t^{k}}{k!}\mid a_{k}\in 
\mathbb{C}
\right\} $, $\mathcal{P}=%
\mathbb{C}
\left[ x\right] $ and $\mathcal{P}^{\ast }$ denotes the vector space of all
linear functional on $\mathcal{P}$, $\left\langle L\mid p\left( x\right)
\right\rangle $ denotes the action of the linear functional $L$ on the
polynomial $p\left( x\right) $, and it is well-known that the vector space
operation on $\mathcal{P}^{\ast }$ is defined by 
\begin{eqnarray*}
\left\langle L+M\mid p\left( x\right) \right\rangle  &=&\left\langle L\mid
p\left( x\right) \right\rangle +\left\langle M\mid p\left( x\right)
\right\rangle , \\
\left\langle cL\mid p\left( x\right) \right\rangle  &=&c\left\langle L\mid
p\left( x\right) \right\rangle ,
\end{eqnarray*}%
where $c$ is some constant in $%
\mathbb{C}
$ (for details, see [10, 11, 14, 17]).

The formal power series are known by the rule%
\begin{equation*}
f\left( t\right) =\sum_{k=0}^{\infty }a_{k}\frac{t^{k}}{k!}\in \mathcal{F}
\end{equation*}%
which defines a linear functional on $\mathcal{P}$ as $\left\langle f\left(
t\right) \mid x^{n}\right\rangle =a_{n}$ for all $n\geq 0$ (for details, 
see [10, 11, 14, 17]). Additionally,%
\begin{equation}
\left\langle t^{k}\mid x^{n}\right\rangle =n!\delta _{n,k},
\label{euation 2}
\end{equation}%
where $\delta _{n,k}$ is the Kronecker symbol. When we take $f_{L}\left(
t\right) =\sum_{k=0}^{\infty }\left\langle L\mid x^{k}\right\rangle \frac{%
t^{k}}{k!},$ then we obtain $\left\langle f_{L}\left( t\right) \mid
x^{n}\right\rangle =\left\langle L\mid x^{n}\right\rangle $ and so as linear
functionals $L=f_{L}\left( t\right) $ (see [10, 11, 14, 17]). Additional,
the map $L\rightarrow f_{L}\left( t\right) $ is a vector space isomorphism
from $\mathcal{P}^{\ast }$ onto $\mathcal{F}$. Henceforth, $\mathcal{F}$
will denote both the algebra of the formal power series in $t$ and the
vector space of all linear functionals on $\mathcal{P}$, and so an element $%
f\left( t\right) $ of $\mathcal{F}$ will be thought of as both a formal
power series and a linear functional. $\mathcal{F}$ will be called as umbral
algebra ( see [10, 11, 14, 17]).

Also, the evaluation functional for $y$ in $%
\mathbb{C}
$ is defined to be power series $e^{yt}$. We can write that $\left\langle
e^{yt}\mid x^{n}\right\rangle =y^{n}$ and so $\left\langle e^{yt}\mid
p\left( x\right) \right\rangle =p\left( y\right) $ ( see [10, 11, 12, 14,
17]). We want to note that for all $f\left( t\right) $ in $\mathcal{F}$%
\begin{equation}
f\left( t\right) =\sum_{k=0}^{\infty }\left\langle f\left( t\right) \mid
x^{k}\right\rangle \frac{t^{k}}{k!}  \label{euation 3}
\end{equation}%
and for all polynomial $p\left( x\right) $, 
\begin{equation}
p\left( x\right) =\sum_{k=0}^{\infty }\left\langle t^{k}\mid p\left(
x\right) \right\rangle \frac{x^{k}}{k!},  \label{euation 4}
\end{equation}%
(for details,  see [10, 11, 14, 17]). The order $o\left( f\left( t\right)
\right) $ of the power series $f\left( t\right) \neq 0$ is the smallest
integer $k$ for which $a_{k}$ does not vanish. It is considered $o\left(
f\left( t\right) \right) =\infty $ if $f\left( t\right) =0$. We see that $%
o\left( f\left( t\right) g\left( t\right) \right) =o\left( f\left( t\right)
\right) +o\left( g\left( t\right) \right) $ and $o\left( f\left( t\right)
+g\left( t\right) \right) \geq \min \left\{ o\left( f\left( t\right) \right)
,o\left( g\left( t\right) \right) \right\} $. The series $f\left( t\right) $
has a multiplicative inverse, denoted by $f\left( t\right) ^{-1}$ or $\frac{1%
}{f\left( t\right) }$, if and only if $o\left( f\left( t\right) \right) =0$.
Such series is called an invertible series. A series $f\left( t\right) $ for
which $o\left( f\left( t\right) \right) =1$ is called a delta series ( see
[10, 11, 12, 14, 17]). For $f\left( t\right) ,g\left( t\right) \in \mathcal{F%
}$, we have $\left\langle f\left( t\right) g\left( t\right) \mid p\left(
x\right) \right\rangle =\left\langle f\left( t\right) \mid g\left( t\right)
p\left( x\right) \right\rangle $.

A delta series $f\left( t\right) $ has a compositional inverse $\overline{f}%
\left( t\right) $ such that $f\left( \overline{f}\left( t\right) \right) =%
\overline{f}\left( f\left( t\right) \right) =t$.

For $f\left( t\right) ,g\left( t\right) \in \mathcal{F}$ , we have $%
\left\langle f\left( t\right) g\left( t\right) \mid p\left( x\right)
\right\rangle =\left\langle f\left( t\right) \mid g\left( t\right) p\left(
x\right) \right\rangle $. By (\ref{euation 3}), we have%
\begin{equation}
p^{\left( k\right) }\left( x\right) =\frac{d^{k}p\left( x\right) }{dx^{k}}%
=\sum_{l=k}^{\infty }\frac{\left\langle t^{l}\mid p\left( x\right)
\right\rangle }{l!}l\left( l-1\right) \cdots \left( l-k+1\right) x^{l-k}%
\text{.}  \label{euation 5}
\end{equation}

Thus, we see that%
\begin{equation}
p^{\left( k\right) }\left( 0\right) =\left\langle t^{k}\mid p\left( x\right)
\right\rangle =\left\langle 1\mid p^{\left( k\right) }\left( x\right)
\right\rangle \text{.}  \label{euation 6}
\end{equation}

By (\ref{euation 5}), we get%
\begin{equation}
t^{k}p\left( x\right) =p^{\left( k\right) }\left( x\right) =\frac{%
d^{k}p\left( x\right) }{dx^{k}}\text{.}  \label{euation 7}
\end{equation}

So, we have%
\begin{equation}
e^{yt}p\left( x\right) =p\left( x+y\right) \text{.}  \label{euation 8}
\end{equation}

Let $S_{n}\left( x\right) $ be a polynomial with $\deg S_{n}\left( x\right)
=n$. Let $f\left( t\right) $ be a delta series and let $g\left( t\right) $
be an invertible series. Then there exists a unique sequence $S_{n}\left(
x\right) $ of polynomials such that $\left\langle g\left( t\right) f\left(
t\right) ^{k}\mid S_{n}\left( x\right) \right\rangle =n!\delta _{n,k}$ for
all $n,k\geq 0$. The sequence $S_{n}\left( x\right) $ is called the sheffer
sequence for $\left( g\left( t\right) ,f\left( t\right) \right) $ or that $%
S_{n}\left( t\right) $ is sheffer for $\left( g\left( t\right) ,f\left(
t\right) \right) $.

The sheffer sequence for $\left( 1,f\left( t\right) \right) $ is called the
associated sequence for $f\left( t\right) $ or $S_{n}\left( x\right) $ is
associated with $f\left( t\right) $. The sheffer sequence for $\left(
g\left( t\right) ,t\right) $ is called the appell sequence for $g\left(
t\right) $ or $S_{n}\left( x\right) $ is Appell for $g\left( t\right) $.

Let $p\left( x\right) \in \mathcal{P}$. Then we have%
\begin{eqnarray}
\left\langle \frac{e^{yt}-1}{t}\mid p\left( x\right) \right\rangle 
&=&\int_{0}^{y}p\left( u\right) du,  \notag \\
\left\langle f\left( t\right) \mid xp\left( x\right) \right\rangle 
&=&\left\langle \partial _{t}f\left( t\right) \mid p\left( x\right)
\right\rangle =\left\langle f%
{\acute{}}%
\left( t\right) \mid p\left( x\right) \right\rangle ,  \label{euation 9} \\
\left\langle e^{yt}-1\mid p\left( x\right) \right\rangle  &=&p\left(
y\right) -p\left( 0\right) ,\text{ ( see [10, 11, 14, 17]).}  \notag
\end{eqnarray}

Let $S_{n}\left( x\right) $ be sheffer for $\left( g\left( t\right) ,f\left(
t\right) \right) $. Then%
\begin{eqnarray}
h\left( t\right) &=&\sum_{k=0}^{\infty }\frac{\left\langle h\left( t\right)
\mid S_{k}\left( x\right) \right\rangle }{k!}g\left( t\right) f\left(
t\right) ^{k},\text{ }h\left( t\right) \in \mathcal{F}  \notag \\
p\left( x\right) &=&\sum_{k=0}^{\infty }\frac{\left\langle g\left( t\right)
f\left( t\right) ^{k}\mid p\left( x\right) \right\rangle }{k!}S_{k}\left(
x\right) ,\text{ }p\left( x\right) \in \mathcal{P},  \notag \\
\frac{1}{g\left( \overline{f}\left( t\right) \right) }e^{y\overline{f}\left(
t\right) } &=&\sum_{k=0}^{\infty }S_{k}\left( y\right) \frac{t^{k}}{k!},%
\text{ for all }y\in 
\mathbb{C}
,  \label{euation 10} \\
f\left( t\right) S_{n}\left( x\right) &=&nS_{n-1}\left( x\right) \text{.} 
\notag
\end{eqnarray}

Let $a_{1},\cdots ,a_{r},b_{1},\cdots ,b_{r}$ be positive integers. Kim and
Rim \cite{Kim1} defined the generating function for multivariate $q$-Euler
polynomials as follows:%
\begin{gather}
F_{q}\left( t,x\mid a_{1},\cdots ,a_{r};b_{1},\cdots ,b_{r}\right)
=\sum_{n=0}^{\infty }E_{n,q}\left( x\mid a_{1},\cdots ,a_{r};b_{1},\cdots
,b_{r}\right) \frac{t^{n}}{n!}  \label{euation 11} \\
=\frac{2^{r}}{\left( q^{b_{1}}e^{a_{1}t}+1\right) \cdots \left(
q^{b_{r}}e^{a_{r}t}+1\right) }e^{xt}\text{.}  \notag
\end{gather}

Note that%
\begin{equation*}
E_{0,q}\left( x\mid a_{1},\cdots ,a_{r};b_{1},\cdots ,b_{r}\right) =\frac{%
2^{r}}{\left[ 2\right] _{q^{b_{1}}}\left[ 2\right] _{q^{b_{2}}}\cdots \left[
2\right] _{q^{b_{r}}}},
\end{equation*}%
where $\left[ x\right] _{q}$ is $q$-extension of $x$ defined by%
\begin{equation*}
\left[ x\right] _{q}=\frac{q^{x}-1}{q-1}=1+q+q^{2}+\cdots +q^{x-1}.
\end{equation*}

We assume that $q\in 
\mathbb{C}
$ with $\left\vert q\right\vert <1$. Also, we want to note that $%
\lim_{q\rightarrow 1}\left[ x\right] _{q}=x$ (see [1-16]). In the special
case, $x=0$, $E_{n,q}\left( 0\mid a_{1},\cdots ,a_{r};b_{1},\cdots
,b_{r}\right) :=E_{n,q}\left( a_{1},\cdots ,a_{r};b_{1},\cdots ,b_{r}\right) 
$ are called multivariate $q$-Euler numbers. By (\ref{euation 11}), we
procure the following:%
\begin{equation}
E_{n,q}\left( x\mid a_{1},\cdots ,a_{r};b_{1},\cdots ,b_{r}\right)
=\sum_{k=0}^{n}\binom{n}{k}x^{k}E_{n-k,q}\left( a_{1},\cdots
,a_{r};b_{1},\cdots ,b_{r}\right) \text{.}  \label{euation 12}
\end{equation}

Kim et al \cite{Kim2} studied some interesting identities for
Frobenius-Euler polynomials arising from umbral calculus. They derived not
only new but also fascianting identities in modern classical umbral
calculus. 

By the same motivation, we also get numerous identities for multivariate $q$%
-Euler polynomials by utilizing from the umbral calculus.

\section{\textbf{On the multivariate }$q$\textbf{-Euler polynomials arising
from umbral calculus}}

Assume that $S_{n}\left( x\right) $ is an appell sequence for $g\left(
t\right) $, by (\ref{euation 10}), we have 
\begin{equation}
\frac{1}{g\left( t\right) }x^{n}=S_{n}\left( x\right) \text{ if and only if }%
x^{n}=g\left( t\right) S_{n}\left( x\right) \text{, }\left( n\geq 0\right) .
\label{euation 13}
\end{equation}

Let us take 
\begin{equation*}
g\left( t\mid a_{1},\cdots ,a_{r};b_{1},\cdots ,b_{r}\right) =\frac{\left(
q^{b_{1}}e^{a_{1}t}+1\right) \cdots \left( q^{b_{r}}e^{a_{r}t}+1\right) }{%
2^{r}}\in \mathcal{F}\text{.}
\end{equation*}

Then we readily see that $g\left( t\mid a_{1},\cdots ,a_{r};b_{1},\cdots
,b_{r}\right) $ is an invertible series. By (\ref{euation 13}), we have 
\begin{equation}
\sum_{n=0}^{\infty }E_{n,q}\left( x\mid a_{1},\cdots ,a_{r};b_{1},\cdots
,b_{r}\right) \frac{t^{n}}{n!}=\frac{1}{g\left( t\mid a_{1},\cdots
,a_{r};b_{1},\cdots ,b_{r}\right) }e^{xt}\text{.}  \label{euation 14}
\end{equation}

By (\ref{euation 14}), we procure the following%
\begin{equation}
\frac{1}{g\left( t\mid a_{1},\cdots ,a_{r};b_{1},\cdots ,b_{r}\right) }%
x^{n}=E_{n,q}\left( x\mid a_{1},\cdots ,a_{r};b_{1},\cdots ,b_{r}\right) 
\text{.}  \label{euation 15}
\end{equation}

Also, by (\ref{euation 10}), we have%
\begin{gather}
tE_{n,q}\left( x\mid a_{1},\cdots ,a_{r};b_{1},\cdots ,b_{r}\right) =E%
{\acute{}}%
_{n,q}\left( x\mid a_{1},\cdots ,a_{r};b_{1},\cdots ,b_{r}\right)
\label{euation 16} \\
=nE_{n-1,q}\left( x\mid a_{1},\cdots ,a_{r};b_{1},\cdots ,b_{r}\right) \text{%
.}  \notag
\end{gather}

By (\ref{euation 15}) and (\ref{euation 16}), we have the following
proposition.

\begin{proposition}
\label{Proposition 1}\textrm{For }$n\geq 0$\textrm{, }$E_{n,q}\left( x\mid
a_{1},\cdots ,a_{r};b_{1},\cdots ,b_{r}\right) $\textrm{\ is an Appell
sequence for} 
\begin{equation*}
g\left( t\mid a_{1},\cdots ,a_{r};b_{1},\cdots ,b_{r}\right) =\frac{\left(
q^{b_{1}}e^{a_{1}t}+1\right) \cdots \left( q^{b_{r}}e^{a_{r}t}+1\right) }{%
2^{r}}\text{.}
\end{equation*}
\end{proposition}

By (\ref{euation 11}), we see that%
\begin{align}
\sum_{n=1}^{\infty }E_{n,q}\left( x\mid a_{1},\cdots ,a_{r};b_{1},\cdots
,b_{r}\right) \frac{t^{n}}{n!}& =\frac{xge^{xt}-g%
{\acute{}}%
e^{xt}}{g^{2}}  \label{euation 17} \\
& =\sum_{n=0}^{\infty }\left( x\frac{1}{g}x^{n}-\frac{g%
{\acute{}}%
}{g}\frac{1}{g}x^{n}\right) \frac{t^{n}}{n!}  \notag
\end{align}%
where we used $g:=g\left( t\mid a_{1},\cdots ,a_{r};b_{1},\cdots
,b_{r}\right) $. Because of (\ref{euation 15}) and (\ref{euation 17}), we
discover the following:%
\begin{gather}
E_{n+1,q}\left( x\mid a_{1},\cdots ,a_{r};b_{1},\cdots ,b_{r}\right) 
\label{euation 18} \\
=xE_{n,q}\left( x\mid a_{1},\cdots ,a_{r};b_{1},\cdots ,b_{r}\right) -\frac{g%
{\acute{}}%
}{g}E_{n,q}\left( x\mid a_{1},\cdots ,a_{r};b_{1},\cdots ,b_{r}\right) \text{%
.}  \notag
\end{gather}

Therefore, we deduce the following theorem.

\begin{theorem}
\label{Theorem 1}\textrm{Let }$g:=g\left( t\mid a_{1},\cdots
,a_{r};b_{1},\cdots ,b_{r}\right) =\frac{\left( q^{b_{1}}e^{a_{1}t}+1\right)
\cdots \left( q^{b_{r}}e^{a_{r}t}+1\right) }{2^{r}}\in F$\textrm{. Then we
have for }$n\geq 0:$%
\begin{equation}
E_{n+1,q}\left( x\mid a_{1},\cdots ,a_{r};b_{1},\cdots ,b_{r}\right) =\left(
x-\frac{g%
{\acute{}}%
}{g}\right) E_{n,q}\left( x\mid a_{1},\cdots ,a_{r};b_{1},\cdots
,b_{r}\right) \text{.}  \label{euation 19}
\end{equation}
\end{theorem}

From (\ref{euation 11}), we derive that%
\begin{gather}
\sum_{n=0}^{\infty }\left( q^{b_{r}}E_{n,q}\left( x+a_{r}\mid a_{1},\cdots
,a_{r};b_{1},\cdots ,b_{r}\right) +E_{n,q}\left( x\mid a_{1},\cdots
,a_{r};b_{1},\cdots ,b_{r}\right) \right) \frac{t^{n}}{n!}
\label{euation 20} \\
=2\sum_{n=0}^{\infty }E_{n,q}\left( x\mid a_{1},\cdots ,a_{r-1};b_{1},\cdots
,b_{r-1}\right) \frac{t^{n}}{n!}\text{.}  \notag
\end{gather}

By comparing the coefficients in the both sides of $\frac{t^{n}}{n!}$ on the
above, we procure the following%
\begin{gather}
2E_{n,q}\left( x\mid a_{1},\cdots ,a_{r-1};b_{1},\cdots ,b_{r-1}\right)
=q^{b_{r}}E_{n,q}\left( x+a_{r}\mid a_{1},\cdots ,a_{r};b_{1},\cdots
,b_{r}\right)   \label{euation 21} \\
+E_{n,q}\left( x\mid a_{1},\cdots ,a_{r};b_{1},\cdots ,b_{r}\right) \text{.}
\notag
\end{gather}

From theorem \ref{Theorem 1}, we get the following equation%
\begin{gather}
gE_{n+1,q}\left( x\mid a_{1},\cdots ,a_{r};b_{1},\cdots ,b_{r}\right) 
\label{euation 22} \\
=gxE_{n,q}\left( x\mid a_{1},\cdots ,a_{r};b_{1},\cdots ,b_{r}\right) -g%
{\acute{}}%
E_{n,q}\left( x\mid a_{1},\cdots ,a_{r};b_{1},\cdots ,b_{r}\right) .  \notag
\end{gather}

By using (\ref{euation 21}) and (\ref{euation 22}), we obtain the following
theorem.

\begin{theorem}
\textrm{For }$n\geq 0$\textrm{, then we have}%
\begin{gather}
2E_{n,q}\left( x\mid a_{1},\cdots ,a_{r-1};b_{1},\cdots ,b_{r-1}\right)
=q^{b_{r}}E_{n,q}\left( x+a_{r}\mid a_{1},\cdots ,a_{r};b_{1},\cdots
,b_{r}\right)   \label{euation 23} \\
+E_{n,q}\left( x\mid a_{1},\cdots ,a_{r};b_{1},\cdots ,b_{r}\right) \text{.}
\notag
\end{gather}
\end{theorem}

Now, we consider that%
\begin{eqnarray*}
&&\int_{x}^{x+y}E_{n,q}\left( u\mid a_{1},\cdots ,a_{r};b_{1},\cdots
,b_{r}\right) du \\
&=&\frac{1}{n+1}\left( E_{n,q}\left( x+y\mid a_{1},\cdots
,a_{r};b_{1},\cdots ,b_{r}\right) -E_{n,q}\left( x\mid a_{1},\cdots
,a_{r};b_{1},\cdots ,b_{r}\right) \right) \\
&=&\frac{1}{n+1}\sum_{j=1}^{\infty }\binom{n+1}{j}E_{n+1-j,q}\left( x\mid
a_{1},\cdots ,a_{r};b_{1},\cdots ,b_{r}\right) y^{j} \\
&=&\sum_{j=1}^{\infty }\frac{n\left( n-1\right) \left( n-2\right) \cdots
\left( n-j+2\right) }{j!}E_{n+1-j,q}\left( x\mid a_{1},\cdots
,a_{r};b_{1},\cdots ,b_{r}\right) y^{j} \\
&=&\frac{1}{t}\left( \sum_{j=0}^{\infty }\frac{y^{j}t^{j}}{j!}-1\right)
E_{n,q}\left( x\mid a_{1},\cdots ,a_{r};b_{1},\cdots ,b_{r}\right) \\
&=&\frac{e^{yt}-1}{t}E_{n,q}\left( x\mid a_{1},\cdots ,a_{r};b_{1},\cdots
,b_{r}\right) \text{.}
\end{eqnarray*}

Therefore, we discover the following theorem:

\begin{theorem}
\textrm{For }$n\geq 0$\textrm{, then we have}%
\begin{equation}
\int_{x}^{x+y}E_{n,q}\left( u\mid a_{1},\cdots ,a_{r};b_{1},\cdots
,b_{r}\right) du=\frac{e^{yt}-1}{t}E_{n,q}\left( x\mid a_{1},\cdots
,a_{r};b_{1},\cdots ,b_{r}\right) \text{.}  \label{euation 24}
\end{equation}
\end{theorem}

By (\ref{euation 16}) and proposition 1, we have%
\begin{equation}
t\left\{ \frac{1}{n+1}E_{n+1,q}\left( x\mid a_{1},\cdots ,a_{r};b_{1},\cdots
,b_{r}\right) \right\} =E_{n,q}\left( x\mid a_{1},\cdots ,a_{r};b_{1},\cdots
,b_{r}\right) \text{.}  \label{euation 25}
\end{equation}

Thanks to (\ref{euation 9}), we readily derive the following:%
\begin{eqnarray}
&&\left\langle e^{yt}-1\mid \frac{E_{n+1,q}\left( x\mid a_{1},\cdots
,a_{r};b_{1},\cdots ,b_{r}\right) }{n+1}\right\rangle  \label{euation 26} \\
&=&\left\langle \frac{e^{yt}-1}{t}\mid t\left\{ \frac{E_{n+1,q}\left( x\mid
a_{1},\cdots ,a_{r};b_{1},\cdots ,b_{r}\right) }{n+1}\right\} \right\rangle 
\notag \\
&=&\left\langle \frac{e^{yt}-1}{t}\mid E_{n,q}\left( x\mid a_{1},\cdots
,a_{r};b_{1},\cdots ,b_{r}\right) \right\rangle .  \notag
\end{eqnarray}

On account of (\ref{euation 25}) and (\ref{euation 26}), we get%
\begin{gather*}
\left\langle \frac{e^{yt}-1}{t}\mid E_{n,q}\left( x\mid a_{1},\cdots
,a_{r};b_{1},\cdots ,b_{r}\right) \right\rangle =\left\langle e^{yt}-1\mid 
\frac{E_{n+1,q}\left( x\mid a_{1},\cdots ,a_{r};b_{1},\cdots ,b_{r}\right) }{%
n+1}\right\rangle \\
=\frac{1}{n+1}\left\{ E_{n+1,q}\left( y\mid a_{1},\cdots ,a_{r};b_{1},\cdots
,b_{r}\right) -E_{n+1,q}\left( a_{1},\cdots ,a_{r};b_{1},\cdots
,b_{r}\right) \right\} \\
=\int_{0}^{y}E_{n,q}\left( u\mid a_{1},\cdots ,a_{r};b_{1},\cdots
,b_{r}\right) du.
\end{gather*}

Consequently, we obtain the following theorem.

\begin{theorem}
\textrm{For }$n\geq 0$\textrm{, then we have}%
\begin{equation}
\left\langle \frac{e^{yt}-1}{t}\mid E_{n,q}\left( x\mid a_{1},\cdots
,a_{r};b_{1},\cdots ,b_{r}\right) \right\rangle =\int_{0}^{y}E_{n,q}\left(
u\mid a_{1},\cdots ,a_{r};b_{1},\cdots ,b_{r}\right) du\text{.}
\label{euation 27}
\end{equation}
\end{theorem}

Assume that 
\begin{equation*}
\mathcal{P}\left( q\mid a_{1},\cdots ,a_{r};b_{1},\cdots ,b_{r}\right)
=\left\{ p\left( x\right) \in Q\left( q\mid a_{1},\cdots ,a_{r};b_{1},\cdots
,b_{r}\right) \left[ x\right] \mid \deg p\left( x\right) \leq n\right\}
\end{equation*}
is a vector space over $Q\left( q\mid a_{1},\cdots ,a_{r};b_{1},\cdots
,b_{r}\right) $.

For $p\left( x\right) \in \mathcal{P}\left( q\mid a_{1},\cdots
,a_{r};b_{1},\cdots ,b_{r}\right) $, let us consider 
\begin{equation}
p\left( x\right) =\sum_{k=0}^{n}b_{k}E_{k,q}\left( x\mid a_{1},\cdots
,a_{r};b_{1},\cdots ,b_{r}\right) \text{.}  \label{euation 28}
\end{equation}

By proposition \ref{Proposition 1}, $E_{n,q}\left( u\mid a_{1},\cdots
,a_{r};b_{1},\cdots ,b_{r}\right) $ is an appell sequence for 
\begin{equation*}
g:=g\left( t\mid a_{1},\cdots ,a_{r};b_{1},\cdots ,b_{r}\right) =\frac{%
\left( q^{b_{1}}e^{a_{1}t}+1\right) \cdots \left(
q^{b_{r}}e^{a_{r}t}+1\right) }{2^{r}}\text{.}
\end{equation*}

Thus we have%
\begin{equation}
\left\langle g\left( t\mid a_{1},\cdots ,a_{r};b_{1},\cdots ,b_{r}\right)
t^{k}\mid E_{n,q}\left( x\mid a_{1},\cdots ,a_{r};b_{1},\cdots ,b_{r}\right)
\right\rangle =n!\delta _{n,k}\text{.}  \label{euation 29}
\end{equation}

From (\ref{euation 28}) and (\ref{euation 29}), we compute%
\begin{gather}
\left\langle g\left( t\mid a_{1},\cdots ,a_{r};b_{1},\cdots ,b_{r}\right)
t^{k}\mid p\left( x\right) \right\rangle =\sum_{l=0}^{n}b_{l}\left\langle
gt^{k}\mid E_{l,q}\left( x\mid a_{1},\cdots ,a_{r};b_{1},\cdots
,b_{r}\right) \right\rangle  \label{euation 30} \\
=\sum_{l=0}^{n}b_{l}l!\delta _{l,k}=k!b_{k}\text{.}  \notag
\end{gather}

Thus, by (\ref{euation 30}), we derive%
\begin{eqnarray}
b_{k} &=&\frac{1}{k!}\left\langle gt^{k}\mid p\left( x\right) \right\rangle
\label{euation 31} \\
&=&\frac{1}{2^{r}k!}\left\langle \left( q^{b_{1}}e^{a_{1}t}+1\right) \cdots
\left( q^{b_{r}}e^{a_{r}t}+1\right) \mid p^{\left( k\right) }\left( x\right)
\right\rangle \text{.}  \notag
\end{eqnarray}

It is not difficult to show the following%
\begin{equation}
\left( q^{b_{1}}e^{a_{1}t}+1\right) \cdots \left(
q^{b_{r}}e^{a_{r}t}+1\right) =\sum_{\underset{k_{1}+k_{2}+\cdots +k_{r}=1}{%
k_{1},\cdots ,k_{r}\geq 0}}q^{\sum_{l=1}^{r}b_{l}k_{l}}e^{t%
\sum_{j=1}^{r}a_{j}k_{j}}\text{.}  \label{euation 32}
\end{equation}

Via the (\ref{euation 31}) and (\ref{euation 32}), we easily see that%
\begin{eqnarray*}
b_{k} &=&\frac{1}{2^{r}k!}\sum_{\underset{k_{1}+k_{2}+\cdots +k_{r}=1}{%
k_{1},\cdots ,k_{r}\geq 0}}q^{\sum_{l=1}^{r}b_{l}k_{l}}\left\langle
e^{t\sum_{j=1}^{r}a_{j}k_{j}}\mid p^{\left( k\right) }\left( x\right)
\right\rangle \\
&=&\frac{1}{2^{r}k!}\sum_{\underset{k_{1}+k_{2}+\cdots +k_{r}=1}{%
k_{1},\cdots ,k_{r}\geq 0}}q^{\sum_{l=1}^{r}b_{l}k_{l}}p^{\left( k\right)
}\left( \sum_{j=1}^{r}a_{j}k_{j}\right) \text{.}
\end{eqnarray*}

As a result, we state the following theorem.

\begin{theorem}
\textrm{For }$p\left( x\right) \in P\left( q\mid a_{1},\cdots
,a_{r};b_{1},\cdots ,b_{r}\right) $\textrm{, when we consider }%
\begin{equation*}
p\left( x\right) =\sum_{k=0}^{n}b_{k}E_{k,q}\left( x\mid a_{1},\cdots
,a_{r};b_{1},\cdots ,b_{r}\right) \text{,}
\end{equation*}%
\textrm{then we have}%
\begin{equation*}
b_{k}=\frac{1}{2^{r}k!}\sum_{\underset{k_{1}+k_{2}+\cdots +k_{r}=1}{%
k_{1},\cdots ,k_{r}\geq 0}}q^{\sum_{l=1}^{r}b_{l}k_{l}}p^{\left( k\right)
}\left( \sum_{j=1}^{r}a_{j}k_{j}\right) \text{,}
\end{equation*}%
\textrm{where }$p^{\left( k\right) }\left( \sum_{j=1}^{r}a_{j}k_{j}\right) =%
\frac{d^{k}p\left( x\right) }{dx^{k}}\mid _{x=\sum_{j=1}^{r}a_{j}k_{j}}$%
\textrm{.}
\end{theorem}

\end{document}